\documentclass[a4paper,onecolumn]{amsart}

	\usepackage{amsmath}
	\usepackage{amsfonts}
	\usepackage{amssymb}
	\usepackage[ansinew]{inputenc}
	\usepackage[usenames,dvipsnames]{pstricks}
	\usepackage{subfigure}
	\usepackage{epsfig}
	\usepackage{pgfplots}
	\usepackage{pst-grad} 
	\usepackage{pst-plot} 
	\usepackage[hyperindex]{hyperref}
	\usepackage{cleveref}
	\usepackage{scrextend}
	\usepackage{amsthm}
	\usepackage{soul}

	\usetikzlibrary{decorations.markings}
	
\changefontsizes[15pt]{11pt}

	\newtheorem{theorem}{Theorem}

\allowdisplaybreaks

    \usepackage[
    top    = 2.00cm,
    bottom = 2.00cm,
    left   = 1.75cm,
    right  = 1.75cm]{geometry}
    
\pgfplotsset{compat=1.13} 

\DeclareFontFamily{U}{mathx}{\hyphenchar\font45}
\DeclareFontShape{U}{mathx}{m}{n}{<->mathx10}{}
\DeclareFontSubstitution{U}{mathx}{m}{n}
\DeclareSymbolFont{mathx}{U}{mathx}{m}{n}
\DeclareMathSymbol{\intop}  {\mathop}{mathx}{"B3} 

\begin{document}

\title{Log-trigonometric integrals and elliptic functions}

\author{Martin Nicholson}

\begin{abstract} A class of log-trigonometric integrals are evaluated in terms of elliptic functions. From this, by using the elliptic integral singular values, one can obtain closed form evaluations of integrals such as
\[
\int\limits_0^{{\pi}/{2}}\ln\left(\cosh\frac{x}{\sqrt{3}}+\cos\frac{\ln \left(2\cos x\right)}{\sqrt{3}}\right)dx=\frac{\pi^2}{8\sqrt{3}}-\frac{\pi}{4}\ln\left(1+\sqrt{3}\right)+\frac{13\pi}{24}\ln 2.
\]
\end{abstract}
\maketitle
\thispagestyle{empty}
\section{Introduction}

The following integrals were calculated in \cite{v}
\begin{align}
    &\int\limits_0^{{\pi}/{2}}\ln\left[x^2+\ln^2(2e^{-a}\cos
    x)\right]dx=\pi\ln \frac{a}{e^b-1},\label{intlog1}\\
    &\int\limits_0^{{\pi}/{2}}\ln\left[x^2+\ln^2(2e^{-a}\cos
    x)\right]\cos2x~dx=\frac{\pi}{2}\left(1-\frac{1}{a}-e^b+\frac{1}{e^b-1}\right),\label{intlog2}\\
    &    \int\limits_{0}^{{\pi}/{2}}\frac{x\sin 2x}{x^2+\ln^2(2e^{-a}\cos x)}~d
    x=\frac{\pi
    }{4}\left(\frac{1}{a^2}+e^{b}-\frac{e^b}{(e^b-1)^2}\right),\label{intlog3}\\
    &\int\limits_{-{\pi}/{2}}^{{\pi}/{2}}\frac{\left(1+e^{2ix}\right)^\gamma}{ix-a+\ln\left(2\cos x\right)}~d x
    =-\frac{\pi}{a}+\pi\frac{e^{(\gamma+1)a}}{e^a-1}H(\ln 2-a),\label{intlog4}
\end{align}
where $a\in\mathbb{R}$, $b=\min\{ a,\ln2\}$, and $H$ is the unit step function. These are log-trigonometric integrals of the type whose study was initiated in the series of papers \cite{glasser,oloa,amdeberhan,bailey,dixit}. The author of the paper \cite{v} also noted that integral \eqref{intlog1} can be used to obtain integrals that can be evaluated in terms of logarithm of Dedekind eta function. However the resulting integrals contained special functions of complex argument. In this paper we modify this approach to obtain integrals of elementary functions of real argument that are evaluated in terms of infinite products or Lambert series. These infinite products and Lambert series can be expressed in terms of elliptic integrals and allow one to obtain closed form evaluation of certain log-trigonometric integrals at particular values of the parameter.

In the following we will use standard notations from the theory of elliptic functions. Let $k$ be modulus, $k'=\sqrt{1-k^2}$ the complementary modulus, and define the complete elliptic integrals of the first and second kind with modulus $k$
$$
K=K(k)=\int\limits_0^{\pi/2}\frac{d\varphi}{\sqrt{1-k^2\sin^2\varphi}},\quad E=E(k)=\int\limits_0^{\pi/2}\sqrt{1-k^2\sin^2\varphi}~{d\varphi},
$$
respectively. Let $K'=K(k')$ and define $\alpha$ according to $\alpha=K'/K$. These notations will be used throught the paper. So that whenever a parameter $\alpha$ is encountered in a formula it is assumed that the corresponding values of $k$, $k'$, $K$, $K'$, $E$ contained in the same formula are determined from the formulas above.

The outline of the paper is as follows. In the next section \ref{main} we formulate and prove main theorems of this paper. These theorems are derived from the integrals (\ref{intlog1}-\ref{intlog4}). However these integrals can be generalized, as shown in Appendix \ref{a} for the case $\gamma=0$ of integral (\ref{intlog4}). Section \ref{other} is devoted to the study of these integrals that are consequences of the integral in Appendix \ref{a}. All integrals in sections \ref{main} and \ref{other} are expressed in terms of elliptic integrals. This means that elliptic integral singular values provide closed form evaluations of log-trigonometric integrals derived in sections \ref{main}, \ref{other}  for particular values of the parameter $\alpha$. Corresponding illustrative examples are considered in section \ref{ex}. Last section \ref{disc} is devoted to discussion of the results presented in this paper. 

\section{Main theorems}\label{main}

\begin{theorem} If $k$, $k'$, and $K$ are defined in terms of $\alpha$ as described in the introduction, then
\begin{equation}\label{1}
    \int\limits_{0}^{{\pi}/{2}}\ln\left(\cosh \frac{x}{\alpha}-\cos \frac{\ln (2 \cos x)}{\alpha}\right)dx=-\frac{\pi^2\alpha}{12}-\frac{\pi}{6}\ln\frac{16kk'K^3\alpha^6}{\pi^3}, \quad \alpha>\frac{\ln 2}{2\pi};
\end{equation}
\begin{equation}\label{2}
    \int\limits_{0}^{{\pi}/{2}}\ln\left(\cosh \frac{x}{\alpha}+\cos \frac{\ln (2 \cos x)}{\alpha}\right)dx=\frac{\pi^2\alpha}{24}+\frac{\pi}{6}\ln\frac{4\sqrt{k}}{k'}, \quad \alpha>\frac{\ln 2}{\pi}.
\end{equation}
\end{theorem}

\noindent
{\it{Proof.}} For $r>\ln 2$ and $n\in \mathbb{N}$, one has from equation \eqref{intlog1}
\begin{equation*}\label{proof1}
    \int\limits_0^{{\pi}/{2}}\ln\left\{\frac{x^2}{r^2n^2}+\left[\frac{\ln(2\cos
    x)}{r n}-1\right]^2\right\}dx=0,
\end{equation*}
\begin{equation*}\label{proof2}
    \int\limits_0^{{\pi}/{2}}\ln\left\{\frac{x^2}{r^2n^2}+\left[\frac{\ln(2\cos
    x)}{r n}+1\right]^2\right\}dx=-\pi\ln(1-e^{-r n}).
\end{equation*}
Now we take the sum of these equations from $n=1$ to infinity. Using the formula
\begin{align*}
\left\{\frac{x^2}{r^2n^2}+\left[\frac{\ln(2\cos
    x)}{r n}-1\right]^2\right\}\left\{\frac{x^2}{r^2n^2}+\left[\frac{\ln(2\cos
    x)}{r n}+1\right]^2\right\}\\
    =1+2\,\frac{x^2-\log ^2(2 \cos x)}{r^2n^2}+\frac{\left[x^2+\log ^2(2 \cos x)\right]^2}{r^4n^4},
\end{align*}
and
$$
\prod_{n=1}^\infty \left\{1+2\,\frac{x^2-\log ^2(2 \cos x)}{r^2n^2}+\frac{\left[x^2+\log ^2(2 \cos x)\right]^2}{r^4n^4}\right\}=\frac{r^2}{2\pi^2}\frac{\cosh \frac{2 \pi  x}{r}-\cos \frac{2 \pi  \ln (2 \cos x)}{r}}{x^2+\log ^2(2 \cos x)},
$$
yields
\begin{equation*}
    \int\limits_{0}^{{\pi}/{2}}\ln\left\{\frac{r^2}{2\pi^2}\frac{\cosh \frac{2 \pi  x}{r}-\cos \frac{2 \pi  \ln (2 \cos x)}{r}}{x^2+\log ^2(2 \cos x)}\right\}dx=-\pi\ln \prod_{n=0}^\infty(1-e^{-r n})\,.
\end{equation*}
Interchanging the order of infinite product and integral can be justified using Fubini's theorem. Using 
\begin{equation*}\label{int3}
    \int\limits_0^{{\pi}/{2}}\ln\left[x^2+\ln^2(2\cos x)\right]dx=0\,,
\end{equation*}
and replacing $r$ by $2\pi\alpha$ one gets
\begin{equation}\label{proof3}
    \int\limits_{0}^{{\pi}/{2}}\ln\left(\cosh \frac{x}{\alpha}-\cos \frac{\ln (2 \cos x)}{\alpha}\right)dx=-\frac{\pi}{2} \ln({2\alpha^2})-\pi\ln \prod_{n=1}^\infty(1-e^{-2\pi\alpha n})\,, \qquad \alpha>\frac{\ln 2}{2\pi};
\end{equation}
\begin{equation}\label{proof4}
    \int\limits_{0}^{{\pi}/{2}}\ln\left(\cosh \frac{x}{\alpha}+\cos \frac{\ln (2 \cos x)}{\alpha}\right)dx=\frac{\pi}{2}\ln2+\pi\ln \prod_{n=1}^\infty(1+e^{-\pi\alpha n})\,, \qquad \alpha>\frac{\ln 2}{\pi}.
\end{equation}
Equation \eqref{proof4} is easily deduced from \eqref{proof3} by taking its linear combinations.
The infinite products in these formulas are calculated in \cite{whittaker}, ch.21, ex.10:
\begin{equation*}\label{e1}
    \prod_{n=1}^\infty(1-e^{-2\pi\alpha n})=e^\frac{\pi\alpha}{12}\left(\frac{2kk'K^3}{\pi^3}\right)^{1/6},
\end{equation*}
\begin{equation*}\label{e2}
    \prod_{n=1}^\infty(1+e^{-\pi\alpha n})=e^\frac{\pi\alpha}{24}\left(\frac{\sqrt{k}}{2k'}\right)^{1/6}.\phantom{.......}
\end{equation*}
This completes the proof.\qed

\begin{theorem}
\begin{equation}\label{3}
    \int\limits_{0}^{{\pi}/{2}}\frac{\cosh \frac{x}{2\alpha}\cos \frac{\ln (2 \cos x)}{2\alpha}}{\cosh \frac{x}{\alpha}+\cos \frac{\ln (2 \cos x)}{\alpha}}\,dx=\frac{\pi}{8}(\alpha+2)-\frac{\alpha}{4}K,\quad \alpha>\frac{\ln 2}{\pi}.
\end{equation}
\end{theorem}

\noindent
{\it{Proof.}} The case $\gamma=0$ of \eqref{intlog4} gives for $r>\ln 2$
\begin{equation*}
    \int\limits_{-{\pi}/{2}}^{{\pi}/{2}}\frac{1}{ix+r(2n+1)+\ln\left(2\cos x\right)}~d x
    =\frac{\pi}{r(2n+1)}, \quad n=-1,-2,-3,...
\end{equation*}
\begin{equation*}
    \int\limits_{-{\pi}/{2}}^{{\pi}/{2}}\frac{1}{ix+r(2n+1)+\ln\left(2\cos x\right)}~d x
    =\frac{\pi}{r(2n+1)}-\pi\frac{1}{e^{r(2n+1)}-1}, \quad n=0,1,2,...
\end{equation*}
Summing these equations with the help of
\begin{equation*}
\sum _{n=-\infty }^{\infty } \frac{(-1)^n}{i x+r(2n+1)+\log (2 \cos x)}=\frac{\pi}{r}\frac{\cosh \frac{\pi  x}{2r}\cos \frac{ \pi  \ln (2 \cos x)}{2r}}{\cosh \frac{\pi  x}{r}+\cos \frac{\pi  \ln (2 \cos x)}{r}},
\end{equation*}
and replacing $r$ by $\pi \alpha$, one is immediately lead to
\begin{equation*}
    \int\limits_{0}^{{\pi}/{2}}\frac{\cosh \frac{x}{2\alpha}\cos \frac{\ln (2 \cos x)}{2\alpha}}{\cosh \frac{x}{\alpha}+\cos \frac{\ln (2 \cos x)}{\alpha}}\,dx=\frac{\pi}{4}-\frac{\pi\alpha}{2} \sum _{n=0}^{\infty} \frac{(-1)^n}{e^{\pi\alpha(2 n+1)}-1}\,,\qquad \alpha>\frac{\ln 2}{\pi}.
\end{equation*}
The interchanging of summation and integration is easily justified by Fubini's theorem. The infinite series in this expression is calculated in \cite{whittaker}, ch.22.6:
\[
\pushQED{\qed} 
    \sum _{n=0}^{\infty} \frac{(-1)^n}{e^{\pi\alpha(2 n+1)}-1}=\frac{K}{2\pi}-\frac{1}{4}.\qedhere
\popQED
\] 

\begin{theorem}
\begin{equation}\label{4}
    \int\limits_{0}^{{\pi}/{2}}\frac{\sin 2x\,\sinh \frac{x}{\alpha}}{\cosh \frac{x}{\alpha}-\cos \frac{\ln (2 \cos x)}{\alpha}}\,dx=\frac{13\pi\alpha}{48}+\frac{\pi}{24\alpha}+\frac{\pi\alpha}{4\tanh 2\pi\alpha}+
\frac{\alpha}{4\pi}\left(E-\frac{2-k^2}{3}K\right)K, \quad \alpha>\frac{\ln 2}{2\pi};
\end{equation}
\begin{equation}\label{5}
    \int\limits_{0}^{{\pi}/{2}}\frac{\sin 2x\,\sinh \frac{x}{\alpha}}{\cosh \frac{x}{\alpha}+\cos \frac{\ln (2 \cos x)}{\alpha}}\,dx=\frac{\pi}{8\alpha}+\frac{\pi\alpha}{4\sinh \pi\alpha}+
\frac{\alpha}{4\pi}\left(E-K\right)K, \quad \alpha>\frac{\ln 2}{\pi}.
\end{equation}
\end{theorem}

\noindent
{\it{Proof.}} It is convenient to rewrite formula \eqref{intlog3} as
\begin{equation}\label{sine}
    \int\limits_{-{\pi}/{2}}^{{\pi}/{2}}\frac{\sin 2x}{ix+\ln\left(2\cos x\right)-a}\,d x
    =\frac{\pi i}{2}\left(\frac{1}{a^2}+e^{b}-\frac{e^b}{(e^b-1)^2}\right).
\end{equation}
Thus for $r>\ln 2$ and $n\in\mathbb{N}$
$$
\int\limits_{-{\pi}/{2}}^{{\pi}/{2}}\frac{\sin 2x}{ix+\ln\left(2\cos x\right)-rn}\,d x
    =\frac{\pi i}{2r^2n^2},
$$
$$
\int\limits_{-{\pi}/{2}}^{{\pi}/{2}}\frac{\sin 2x}{ix+\ln\left(2\cos x\right)+rn}\,d x
    =\frac{\pi i}{2}\left(\frac{1}{r^2n^2}+e^{-rn}-\frac{e^{-rn}}{(1-e^{-rn})^2}\right).
$$
From this it follows with the help of the summation
\begin{equation*}
    \frac{1}{i x+\log (2 \cos x)}+2\sum _{n=1}^{\infty } \frac{i x+\log (2 \cos x)}{[i x+\log (2 \cos x)]^2-r^2n^2}=\frac{\pi}{r}\frac{\sin\frac{2\pi  \ln (2 \cos x)}{r}-i\sinh \frac{2\pi  x}{r}}{\cosh \frac{2\pi  x}{r}-\cos \frac{2\pi  \ln (2 \cos x)}{r}}
\end{equation*}
and integration (\cite{v}, eq. (47))
\begin{equation}\label{sine0}
    \int\limits_{-{\pi}/{2}}^{{\pi}/{2}}\frac{\sin 2x}{ix+\ln\left(2\cos x\right)}\,d x
    =\frac{13\pi}{24},
\end{equation}
that
\begin{equation*}
        \int\limits_{0}^{{\pi}/{2}}\frac{\sin 2x\,\sinh \frac{2\pi  x}{r}}{\cosh \frac{2\pi x}{r}-\cos \frac{2\pi  \ln (2 \cos x)}{r}}\,dx=\frac{\pi^2}{12r }+\frac{13r}{48}+\frac{r}{4(e^{r}-1)}-\frac{r}{16} \sum _{n=0}^{\infty} \frac{1}{\sinh^2{\frac{r n}{2}}}.
\end{equation*}
After change of parameter $r=2\pi\alpha$, this identity takes the form
\begin{equation*}\label{x}
    \int\limits_{0}^{{\pi}/{2}}\frac{\sin 2x\,\sinh \frac{x}{\alpha}}{\cosh \frac{x}{\alpha}-\cos \frac{\ln (2 \cos x)}{\alpha}}\,dx=\frac{\pi}{24\alpha }+\frac{13\pi\alpha}{24}+\frac{\pi\alpha}{2(e^{2\pi\alpha}-1)}-\frac{\pi\alpha}{8} \sum _{n=1}^{\infty} \frac{1}{\sinh^2{\pi\alpha n}}\,, \qquad \alpha>\frac{\ln 2}{2\pi}.
\end{equation*}
One can easily deduce taking linear combinations of the previous identity that
\begin{equation*}\label{y}
    \int\limits_{0}^{{\pi}/{2}}\frac{\sin 2x\,\sinh \frac{x}{\alpha}}{\cosh \frac{x}{\alpha}+\cos \frac{\ln (2 \cos x)}{\alpha}}\,dx=\frac{\pi}{8\alpha }+\frac{\pi\alpha}{4\sinh\pi\alpha}-\frac{\pi\alpha}{8} \sum _{n=0}^{\infty} \frac{1}{\sinh^2\frac{\pi\alpha (2n+1)}{2}}\,, \qquad \alpha>\frac{\ln 2}{\pi}.
\end{equation*}
The infinite series in these formulas were calculated in \cite{ling} in terms of elliptic integrals
\begin{equation}\label{e4}
    \phantom{...........}\sum _{n=1}^{\infty} \frac{1}{\sinh^2{\pi\alpha n}}=\frac{1}{6}-\frac{2}{\pi^2}KE+\frac{2(2-k^2)}{3\pi^2}K^2,
\end{equation}
\begin{equation}\label{e5}
    \sum _{n=0}^{\infty} \frac{1}{\sinh^2\frac{\pi\alpha (2n+1)}{2}}=\frac{2}{\pi^2}K(K-E).\phantom{....}
\end{equation}
The proof of (\ref{4}) and (\ref{5}) is complete. \qed

\begin{theorem}
\begin{equation}\label{6}
    \int\limits_{0}^{{\pi}/{2}}\frac{\cos 2x\,\sin \frac{\ln (2 \cos x)}{\alpha}}{\cosh \frac{x}{\alpha}-\cos \frac{\ln (2 \cos x)}{\alpha}}\,dx=\frac{11\pi\alpha}{48}-\frac{\pi}{24\alpha}+\frac{\pi\alpha}{4\tanh 2\pi\alpha}-
\frac{\alpha}{4\pi}\left(E-\frac{2-k^2}{3}K\right)K, \quad \alpha>\frac{\ln 2}{2\pi};
\end{equation}
\begin{equation}\label{7}
    \int\limits_{0}^{{\pi}/{2}}\frac{\cos 2x\,\sin \frac{\ln (2 \cos x)}{\alpha}}{\cosh \frac{x}{\alpha}+\cos \frac{\ln (2 \cos x)}{\alpha}}\,dx=\frac{\pi}{8\alpha}-\frac{\pi\alpha}{4\sinh \pi\alpha}+
\frac{\alpha}{4\pi}\left(E-K\right)K, \quad \alpha>\frac{\ln 2}{\pi}.
\end{equation}
\end{theorem}

\noindent
{\it{Proof.}} These identities are proved in a similar manner to the proof of the previous theorem. The starting log trigonometric integral is
\begin{equation}\label{cos}
    \int\limits_{-\pi/2}^{\pi/2}\frac{\cos 2x}{ix+\ln\left(2\cos x\right)-a}\,d x
    =\pi e^a H(\ln 2-a)-\frac{\pi}{2}\left(\frac{1}{a^2}+e^{b}-\frac{e^b}{(e^b-1)^2}\right),
\end{equation}
where $a\in\mathbb{R}$, $b=\min\{ a,\ln2\}$, which is a consequence of (\ref{intlog3}) and the $\gamma=0$ and $\gamma=1$ cases of (\ref{intlog4}). From this, it follows that
\begin{equation*}
    \int\limits_{0}^{{\pi}/{2}}\frac{\cos 2x\,\sin \frac{2\pi  \ln (2 \cos x)}{r}}{\cosh \frac{2\pi  x}{r}-\cos \frac{2\pi  \ln (2 \cos x)}{r}}\,dx=\frac{11 r}{48}-\frac{\pi ^2}{12r}+\frac{r}{4 \left(e^r-1\right)}+\frac{r}{16}\sum _{n=1}^{\infty} \frac{1}{\sinh^2{\frac{r n}{2}}},\quad r>\ln 2.
\end{equation*}
Thus
\begin{equation}\label{6a}
    \int\limits_{0}^{{\pi}/{2}}\frac{\cos 2x\,\sin \frac{\ln (2 \cos x)}{\alpha}}{\cosh \frac{x}{\alpha}-\cos \frac{\ln (2 \cos x)}{\alpha}}\,dx=\frac{11 \pi\alpha}{24}-\frac{\pi}{24 \alpha}+\frac{\pi\alpha}{2 \left(e^{2\pi\alpha}-1\right)}+\frac{\pi\alpha}{8}\sum _{n=1}^{\infty} \frac{1}{\sinh^2{\pi\alpha n}}\,, \qquad \alpha>\frac{\ln 2}{2\pi};
\end{equation}
\begin{equation}\label{7a}
    \int\limits_{0}^{{\pi}/{2}}\frac{\cos 2x\,\sin \frac{\ln (2 \cos x)}{\alpha}}{\cosh \frac{x}{\alpha}+\cos \frac{\ln (2 \cos x)}{\alpha}}\,dx=\frac{\pi}{8\alpha }-\frac{\pi\alpha}{4\sinh\pi\alpha}-\frac{\pi\alpha}{8} \sum _{n=0}^{\infty} \frac{1}{\sinh^2\frac{\pi\alpha (2n+1)}{2}}\,, \qquad \alpha>\frac{\ln 2}{\pi}.
\end{equation}
(\ref{7a}) is a consequence of (\ref{6a}). The infinite series in these formulas are the same as in the previous theorem, eqs. (\ref{e4}) and (\ref{e5}). \qed

\section{Some other integrals}\label{other}

\begin{theorem}
\begin{equation}
\int\limits_0^{{\pi}}\frac{ \sinh\frac{ 4x-\pi}{\alpha}}{\cosh{\frac{ 4x-\pi}{\alpha}}-\cos\frac{4\ln (2\sin x)}{\alpha}}\,dx=\pi  \coth \frac{\pi }{2\alpha}-\frac{\pi \alpha}{8(\sqrt{2} -1)}- \frac{\alpha K}{4}\big(1+\sqrt{2+2k}\big),\quad \alpha>\frac{\ln 2}{\pi};
\end{equation}
\begin{equation}
\int\limits_0^{{\pi}}\frac{ \sinh\frac{ 4x-\pi}{\alpha}}{\cosh{\frac{ 4x-\pi}{\alpha}}+\cos\frac{4\ln (2\sin x)}{\alpha}}\,dx=\pi  \tanh \frac{\pi }{2\alpha}- \frac{\alpha kK}{4}\big(1+\sqrt{2+2/k}\big),\quad \alpha>\frac{\ln 4}{\pi}.
\end{equation}
\end{theorem}

\noindent
{\it{Proof.}}
The case $\theta=\pi/4$ of equation (\ref{theta}) is
\begin{equation*}\label{}
    \int\limits_{0}^{{\pi}}\frac{dx}{i(x-\pi/4)-a+\ln\left(2\sin x\right)}
    =\frac{4\pi}{{\pi i}-4a}+\frac{\pi}{1-e^{\pi i/4-a}}H(\ln 2-2a),\quad a\in\mathbb{R}.
\end{equation*}
Replacing $a$ with $\pi\alpha n/2$, where $\alpha>\frac{\ln 2}{\pi}$ and $n\in\mathbb{Z}$, and summing all these equations one can get
\begin{equation*}
\int\limits_0^{{\pi}}\frac{ \sinh\frac{ 4x-\pi}{\alpha}}{\cosh{\frac{ 4x-\pi}{\alpha}}-\cos\frac{4\ln (2\sin x)}{\alpha}}\,dx=\pi  \coth \frac{\pi }{2\alpha}-\frac{\pi \alpha}{4(\sqrt{2} -1)}- \frac{\pi\alpha}{4}\sum _{n=1}^{\infty} \frac{1}{ \sqrt{2} \cosh \frac{\pi\alpha n}{2}-1},
\end{equation*}
and as a consequence
\begin{equation*}
\int\limits_0^{{\pi}}\frac{ \sinh\frac{ 4x-\pi}{\alpha}}{\cosh{\frac{ 4x-\pi}{\alpha}}+\cos\frac{4\ln (2\sin x)}{\alpha}}\,dx=\pi  \tanh \frac{\pi }{2\alpha}- \frac{\pi\alpha}{4}\sum _{n=0}^{\infty} \frac{1}{ \sqrt{2} \cosh \frac{\pi\alpha (2n+1)}{4}-1},\quad \alpha>\frac{\ln 4}{\pi}.
\end{equation*}
The series in these equations is calculated in terms of elliptic integrals in the appendix \ref{b}. \qed

Another case is $\theta=\frac{\pi}{3}$. It is interesting in that it has a weaker restriction on the parameter $\alpha$ than other formulas obtained above. In terms of the base of the elliptic functions this restriction is $q<1$:

\begin{theorem} If $\alpha>0$, then
\begin{equation}
\int\limits_{0}^{{\pi}}\frac{ \sinh{\frac{\pi-6x}{2\alpha}}}{\cosh{\frac{\pi-6x}{2\alpha}}+\cos\frac{3\ln (2\sin x)}{\alpha}}\,dx=\frac{\alpha kK}{\sqrt{3}}\,\mathrm{cn}\left(\tfrac{iK'}{3},k\right)-\pi  \tanh \frac{\pi }{2\alpha}.
\end{equation}
\end{theorem}

\noindent
{\it{Proof.}} The proof is similar to the proof of the previous theorem. We get 
\begin{equation*}
\int\limits_{-\pi/2}^{{\pi}/2}\frac{ \sinh{\left(\frac{ 3x}{\alpha}+\frac{\pi}{\alpha}\right)}\,dx}{\cosh{\left(\frac{ 3x}{\alpha}+\frac{\pi}{\alpha}\right)}+\cos\frac{3\ln (2\cos x)}{\alpha}}=\pi  \tanh \frac{\pi }{2\alpha}-\frac{\pi\alpha}{\sqrt{3}}\sum_{n=0}^{\infty} \frac{1}{2\cosh\frac{\pi\alpha(2n+1)}{3}-1},\quad \alpha>0,
\end{equation*}
and the series in this formula was calculated in the Appendix \ref{b}.\qed

\begin{theorem} If $\alpha>0$, then
\begin{equation}
\int\limits_{0}^{2\pi}\!\arctan\Biggl\{\!\frac{\tanh{\tfrac{\pi-3x}{4\alpha}}}{\cot\!\frac{3\ln (2\sin \frac{x}{2})}{2\alpha}}\!\Biggr\}\cos x\,dx=\frac{\pi  \sqrt{3}}{4 \sinh\frac{\pi\alpha}{3}}-\frac{3 \pi}{2 \alpha}\tanh\frac{\pi }{2 \alpha}+\frac{\sqrt{3}kK}{2}\mathrm{cn}\left(\tfrac{iK'}{3},k\right).
\end{equation}
\end{theorem}

\noindent
{\it{Proof.}} Only an outline of the proof is given here. We start from the generalization of (\ref{cos})
\begin{equation}\label{theta2}
   \int\limits_{-\pi/2}^{{\pi}/2}\frac{\cos 2x\,dx}{i(x+\theta)-a+\ln\left(2\cos x\right)}
    =-\frac{\pi}{2(i\theta - a)^2}+\frac{\pi}{2}\left(e^{a-i\theta}+\frac{e^{a-i\theta}}{\left(1-e^{a-i\theta}\right)^2}\right)H\left[\ln (2\cos\theta)-a\right],
\end{equation}
where $-\frac{\pi}{2}<\theta<\frac{\pi}{2}$, $a\in\mathbb{R}$, and $a\neq \ln (2\cos\theta)$. Its proof is similar to the one considered in Appendix \ref{a} with slight modifications. First we set $\theta=\pi/3$ in this equation, in which case $\ln (2\cos\theta)=0$, then replace $a$ with $(2n+1)r$ where $r>0$, and take the sum from $n=-N-1$ to $n=N$~ $(N\in \mathbb{N})$. It is easily checked that integration and the limit $\lim_{N\to\infty}$ can be interchanged in this case by Fubini's theorem. Hence
\begin{align*}
    \frac{\pi}{2r}\int\limits_{-\pi/2}^{{\pi}/2}\cos 2x\,&\tan\frac{\pi(ix+i\theta+\ln(2 \cos x))}{2 r}\,dx\\
    &=\sum_{n=0}^\infty\frac{\pi}{(i\theta - r(2n+1))^2}-\frac{\pi}{2}\sum_{n=0}^\infty\left(e^{-(2n+1)r-i\theta}+\frac{e^{-(2n+1)r-i\theta}}{\left(1-e^{-(2n+1)r-i\theta}\right)^2}\right).
\end{align*}
On the RHS of this equation, we first integrate wrt to $\theta$ and then take real part, while on the lhs it is easier first to take real part and then to integrate wrt to $\theta$ by using the value of the elementary integral
$$
\int \frac{\sin\frac{\pi  y}{r}}{\cosh\frac{\pi  (\theta+x)}{r}+\cos\frac{\pi  y}{r}} \, d\theta=\frac{2 r}{\pi} \arctan\left\{\tanh\frac{\pi  (\theta+x)}{2 r}\tan\frac{\pi  y}{2 r}\right\}.
$$
The result is
\begin{align*}
    \int\limits_{-\pi/2}^{{\pi}/2}\arctan&\left\{\tanh\frac{\pi  (\pi+3x)}{6 r}\tan\frac{\pi  \ln(2\cos x)}{2 r}\right\}\cos 2x\,dx\\
    &\phantom{phantom}=\frac{\pi^2}{4 r}\tanh\frac{\pi ^2}{6 r}-\frac{\pi \sqrt{3}}{8\sinh r}-\frac{\pi \sqrt{3}}{4}\sum_{n=0}^\infty\frac{1}{2\cosh (2n+1)r-1}.
\end{align*}
The sum on the RHS is calculated in Appendix \ref{b}. \qed

It is clear that a lot of other log-trigonometric evaluations in terms of elliptic integrals can be obtained in this way.

\section{Examples}\label{ex}

Using elliptic integral singular values one can find closed form evaluations of log-trigonometric integrals considered above at certain values of the parameter $\alpha$, e.g.
\begin{align}
    &\int\limits_0^{{\pi}/{2}}\ln\left(\cosh\frac{x}{\sqrt{3}}+\cos\frac{\ln \left(2\cos x\right)}{\sqrt{3}}\right)dx=\frac{\pi^2}{8\sqrt{3}}-\frac{\pi}{4}\ln\left(1+\sqrt{3}\right)+\frac{13\pi}{24}\ln 2,\label{ex2}\\
    &\int\limits_{0}^{{\pi}/{2}}\frac{\cosh \frac{x}{4}\cos \frac{\ln (2 \cos x)}{4}}{\cosh \frac{x}{2}+\cos \frac{\ln (2 \cos x)}{2}}\,dx=\frac{\pi}{2}-\frac{\left(\sqrt{2}+1\right) \Gamma^2\!\left(\frac{1}{4}\right)}{16 \sqrt{2 \pi }},\\
    &\int\limits_{0}^{{\pi}}\frac{ \sinh{\frac{\pi-6x}{2\sqrt{3}}}}{\cosh{\frac{\pi-6x}{2\sqrt{3}}}+\cos\frac{3\ln (2\sin x)}{\sqrt{3}}}\,dx=\frac{\Gamma^3\!\left(\frac{1}{3}\right)}{2^{10/3} \pi }-\pi  \tanh\frac{\pi }{2 \sqrt{3}}.
\end{align}
Equation (\ref{2}) in Theorem $1$ is particularly interesting because in this case $k$ and $k'$ are algebraic. That's why evaluation of the integral (\ref{ex2}) does not contain any gamma functions.

\section{Discussion}\label{disc}
To better understand the integrals calculated in the previous sections, let's consider the integral (\ref{3}) from Theorem $2$ and transform it into another form
\begin{align*}
    I(\alpha)=\int\limits_{0}^{{\pi}/{2}}\frac{\cosh \frac{x}{2\alpha}\cos \frac{\ln (2 \cos x)}{2\alpha}}{\cosh \frac{x}{\alpha}+\cos \frac{\ln (2 \cos x)}{\alpha}}\,dx&=\frac14\int\limits_{-{\pi}/{2}}^{{\pi}/{2}}{\frac{dx}{\cos\left(\frac{ix}{2\alpha}+\frac{\ln (2 \cos x)}{2\alpha}\right)}}\\
    &=\frac14\int\limits_{-{\pi}/{2}}^{{\pi}/{2}}{\frac{dx}{\cos\frac{\ln (1+e^{2ix})}{2\alpha}}}.
\end{align*}
Here, making the substitution $z=\ln (1+e^{2ix})$ as explained in the Appendix \ref{a}, one obtains
\begin{equation}\label{countor_integral}
    I(\alpha)=\frac{1}{8i}\int\limits_C\frac{dz}{(1-e^{-z})\cos\frac{z}{2\alpha}}.
\end{equation}
Contour $C$ is the same as in Appendix \ref{a}. This integral can be evaluated directly by means of residue theorem in terms of a Lambert series.

Thus, what we essentially did in this paper was to rewrite certain contour integrals as integrals of real valued functions over the real interval. 

There are other integrals that can be evaluated in terms of Lambert series, e.g.
\begin{equation}
\int\limits_0^{{\pi}/{2}}\frac{ \sin \frac{\ln (2\cos x)}{\alpha}}{\cosh{\frac{ x}{\alpha}}-\cos\frac{\ln (2\cos x)}{\alpha}}\,dx=\frac{\pi\alpha}{2}-\pi\alpha\sum_{n=1}^\infty\frac{1}{e^{2\pi\alpha n}-1},
\end{equation}
\begin{equation}\label{lambert2}
\int\limits_0^{{\pi}/{2}}\frac{ \sin \frac{\ln (2\cos x)}{\alpha}}{\cosh{\frac{ x}{\alpha}}+\cos\frac{\ln (2\cos x)}{\alpha}}\,dx=\pi\alpha\sum_{n=0}^\infty\frac{1}{e^{\pi\alpha(2n+1)}-1}.
\end{equation}
However these Lambert series can not be expressed in terms of elliptic integrals.

For $\theta=\pi/2$ the discontinuity in the parameter $a$ disappears altogether:
\begin{align*} 
  \int_{0}^{\pi} \frac{x}{x^2+\ln^2(2 e^{a} \sin x)}\,dx 
  & = \, \frac{2 \pi^2}{\pi^2+4a^2},\\ 
  \int_{0}^{\pi} \frac{\ln(2 e^{a} \sin x)}{x^2+\ln^2(2 e^{a} \sin x)}\,dx & =  \frac{4\pi a}{\pi^2+4a^2},\qquad a\in\mathbb{R}.
\end{align*} 
These formulas lead to the analogs of \eqref{lambert2}
\begin{align*} 
  \int_0^{\pi}\frac{\sin\frac{\ln(2\sin x)}{\alpha}}{\cosh\frac{x}{\alpha}+\cos\frac{\ln(2\sin x)}{\alpha}}\,dx&=0,\\ 
  \int_0^{\pi}\frac{\sinh\frac{x}{\alpha}}{\cosh\frac{x}{\alpha}+\cos\frac{\ln(2\sin x)}{\alpha}}\,dx&=\pi\tanh\frac{\pi}{4\alpha},\qquad \alpha\in\mathbb{R},
\end{align*} 
that obviously does not contain any Lambert series.

Another type of integrals can be obtained from log-trigonometric integrals (\ref{intlog1}-\ref{intlog4}) when $|\text{Im}\, a|\ge \frac{\pi}{2}$, e.g.
$$
\int\limits_0^{\pi/2}\frac{\sinh\frac{\ln(2\cos x)}{\alpha}}{\cosh\frac{\ln(2\cos x)}{\alpha}-\cos\frac{x}{\alpha}}=\frac{\pi\alpha}{2},\qquad \alpha\in\mathbb{R}.
$$
These integrals also does not contain any Lambert series.

\appendix

\section{Generalization of equation (\ref{intlog4}) with \texorpdfstring{$\gamma=0$}{Lg}}\label{a}

Here we consider a generalization of (\ref{intlog3}) when $\beta=0$. Namely, it will be proved that
\begin{equation}\label{theta}
    \int\limits_{-\pi/2}^{{\pi}/2}\frac{dx}{i(x+\theta)-a+\ln\left(2\cos x\right)}
    =\frac{\pi}{i\theta - a}+\frac{\pi}{1-e^{i\theta-a}}H\left[\ln (2\cos\theta)-a\right],
\end{equation}
where $-\frac{\pi}{2}<\theta<\frac{\pi}{2}$, $a\in\mathbb{R}$, and $a\neq \ln (2\cos\theta)$.

To do this first rewrite this integral as a contour integral. Let 
\begin{equation}\label{z}
    z=ix+\ln(2\cos x),\qquad -\frac{\pi}{2}<x<\frac{\pi}{2}.
\end{equation}
When $x$ varies from $-\frac{\pi}{2}$ to $\frac{\pi}{2}$ the complex variable $z$ traverses the path $C$ given by the parametric equation $\mathrm{Re}\,z=\ln(2\cos x)$, $\mathrm{Im}\,z=x$. This path is plotted in the figure below (the blue line).
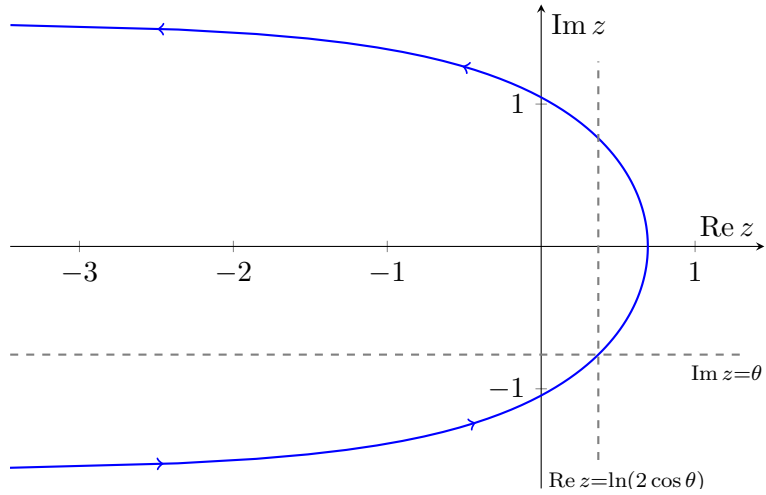
\begin{center}
    \begin{tikzpicture}
\begin{axis}[
    trig format plots=rad,
	axis x line=middle,
	axis y line=middle,
	xtick={-3,-2,-1,0,1},
    ytick={-1,0,1},
	width=11.5cm, height=8cm,
	ymax=1.7, ymin=-1.7, xmax=1.45, xmin=-3.45,
	xlabel=${\mathrm{Re}\,z}$,ylabel=${\mathrm{Im}\,z}$
]
	\addplot [domain=-346.5*pi/700:346.5*pi/700, samples=100, thick, blue, 
	postaction={decorate, decoration={markings,
    mark=at position 0.105 with {\arrow{>};},
    mark=at position 0.32 with {\arrow{>};},
    mark=at position 0.69 with {\arrow{>};},
    mark=at position 0.9 with {\arrow{>};}
      }}
        ] ({ln(2*cos(x))},{(x)});
    \addplot [domain=-3.45:1.3, samples=20, thick, dashed, gray] ({(x)},{-0.76}) node[below,black,pos=0.98] {${\scriptstyle{\text{Im}\,z=\theta}}$};
    \addplot [domain=-1.5:1.3, samples=20, thick, dashed, gray] ({ln(2*cos(0.76))},{(x)}) node[below,black,pos=0.0] {${\scriptstyle{\text{\,\,\,\,\,\,\,\,\,\,\,\,\,\,\,Re}\,z=\ln(2\cos\theta)}}$};
\end{axis}
\end{tikzpicture}
\end{center}
As can be seen, this path extends from $-\infty$ below the line $\text{Im}\, z=0$, passes through the point $(\ln 2,0)$, then extends back to $-\infty$ above the line $\text{Im}\,z=0$. It is easy to solve eq. (\ref{z}) for $x$ in terms of $z$: $x=\frac{1}{2i}\ln(e^z-1)$, with the choice of the branch cut for the complex logarithm as the ray going from $0$ to $-\infty$. Thus
$$
\int\limits_{-\pi/2}^{{\pi}/2}\frac{dx}{i(x+\theta)-a+\ln\left(2\cos x\right)}=\frac{1}{2i}\int\limits_C\frac{dz}{(z+i\theta-a)(1-e^{-z})}.
$$
The integral in the $z$-domain does not have any branching points, therefore there is no need to include any branch cuts in the picture above.

There can be at most two poles of the function $\frac{1}{(z+i\theta-a)(1-e^{-z})}$ inside the contour of integration. The pole $z=0$ is always inside the contour, while the pole $z=a-i\theta$ can be either inside or outside of the contour depending on the values of $a$ and $\theta$. Let's fix the value of $-\frac{\pi}{2}<\theta<\frac{\pi}{2}$ and plot the line $\text{Im}\, z=\theta$ in the picture above (the dashed horizontal line). It has only one point of intersection with path $C$ which is $(\ln(2\cos\theta),\theta)$. Thus if $a<\ln(2\cos\theta)$, then the pole $z=a-i\theta$ is inside the contour $C$. If $a>\ln(2\cos\theta)$ then the pole $z=a-i\theta$ is outside the contour $C$. Now one can easily apply the residue theorem to complete the proof of (\ref{theta}).

\section{Calculating certain sums with hyperbolic functions}\label{b}

Let
$$
S_1(\alpha)=\sum_{n=-\infty}^\infty\frac{1}{\sqrt{2}\cosh\frac{\pi \alpha n}{2}+1},
$$
$$
S_2(\alpha)=\sum_{n=-\infty}^\infty\frac{1}{\sqrt{2}\cosh\frac{\pi \alpha n}{2}-1},
$$
then due to $2\cosh^2x-1=\cosh 2x$ one obtains
$$
S_2(\alpha)-S_1(\alpha)=2\sum_{n=-\infty}^\infty\frac{1}{\cosh \pi \alpha n},
$$
$$
S_2(\alpha)+S_1(\alpha)=2\sqrt{2}\sum_{n=-\infty}^\infty\frac{\cosh\frac{\pi\alpha n}{2}}{\cosh \pi \alpha n}.
$$
Well known formulas from the theory of elliptic functions \cite{whittaker} state that
$$
\sum_{n=-\infty}^\infty\frac{1}{\cosh \pi \alpha n}=\frac{2K}{\pi},\quad\sum_{n=-\infty}^\infty\frac{\cosh\frac{\pi\alpha n}{2}}{\cosh \pi \alpha n}=\frac{2K}{\pi}\sqrt{1+k}.
$$
One can deduce from this by trivial algebra that
$$
S_2(K'/K)=\sum_{n=-\infty}^\infty\frac{1}{\sqrt{2}\cosh\frac{\pi \alpha n}{2}-1}=\frac{2K}{\pi}\big(1+{\sqrt{2+2k}}\big).
$$
Similarly
\begin{equation}
    \sum _{n=0}^{\infty} \frac{1}{ \sqrt{2} \cosh \frac{\pi\alpha (2n+1)}{4}-1}=\frac{kK}{\pi}\big(1+{\sqrt{2+2/k}}\big).
\end{equation}

To calculate the sums we encountered in the proofs of thereoms $6$ and $7$ note that $2\cosh 2x-1=\frac{\cosh 3x}{\cosh x}$. Then
it is easy to deduce that
\begin{equation}
    \sum _{n=0}^{\infty} \frac{1}{2\cosh\frac{\pi\alpha(2n+1)}{3}-1}=\frac{kK }{\pi }\,\text{cn}\left(\tfrac{iK'}{3},k\right).
\end{equation}

\end{document}